# Nonlinear Dynamic Interactions between flow-induced Galloping and Shell-like Buckling

J. Michael T. Thompson[1] & Jan Sieber[2]

[1]*Centre for Applied Dynamics Research, School of Engineering,
University of Aberdeen, AB24 3UE, and
Honorary Fellow, DAMTP, Cambridge.*

[2] College of Engineering, Mathematics and Physical Sciences,
*Harrison Building, University of Exeter, Exeter EX4 4QF, UK.*

## Abstract

For an elastic system that is non-conservative but autonomous, subjected for example to time-independent loading by a steadily flowing fluid (air or water), a dangerous bifurcation, such as a sub-critical bifurcation, or a cyclic fold, will trigger a dynamic jump to one or more remote stable attractors. When there is more than one candidate attractor, the one onto which the structure settles can then be *indeterminate*, being sensitive to infinitesimally small variations in starting conditions or parameters.

   In this paper we develop and study an archetypal model to explore the nonlinear dynamic interactions between galloping at an incipient sub-critical Hopf bifurcation of a structure with shell-like buckling behaviour that is gravity-loaded to approach a sub-critical pitch-fork bifurcation. For the fluid forces, we draw on the aerodynamic coefficients determined experimentally by Novak for the flow around a bluff body of rectangular cross-section. Meanwhile, for the structural component, we consider a variant of the propped-cantilever model that is widely used to illustrate the sub-critical pitchfork: within this model a symmetry-breaking imperfection makes the behaviour generic.

   The compound bifurcation corresponding to simultaneous galloping and buckling is the so-called Takens-Bodganov Cusp. We make a full unfolding of this codimension-3 bifurcation for our archetypal model to explore the adjacent phase-space topologies and their indeterminacies.

## 1. Introduction

The simplest form of pure galloping is exhibited by a bluff body oscillating transversely in a steady wind. With a structural support providing both linear elastic stiffness and linear viscous damping, the theory for this phenomenon was developed by Novak [1] for a series of rectangular cross-sections. Based on experimental fitting to the quasi-static aerodynamic forces, Novak's theory agreed well with his related experimental studies. An excellent modern account of this, and other work, is given in the book by Paidousis *et al* [2]. Note that galloping is essentially a one-mode phenomenon, distinct from flutter which arises in systems with at least two active modes; and even more distinct from vortex resonance which involves a strong interaction with the fluid. Note, though, that in nonlinear dynamics the bifurcations to both galloping and flutter are described as a Hopf bifurcations [3-5].

The essence of Novak's galloping theory was to use the highly non-linear aerodynamic force characteristics obtained by calibration experiments in which a steady wind-stream was directed, at a series of (resultant) angles, towards the stationary rectangular body. The characteristic graph of lateral force versus angle of attack was then approximated by a seventh-order polynomial. Some of Novak's results are summarised in figure 1. Here the lateral force on the rectangular prism, in the direction of the lateral displacement, $x$, due to a wind of velocity, $V$, is $½ \rho a V^2 C_f(\alpha)$ where $\rho$ is the air density, $a$ is the frontal area, and the (small) angle $\alpha$ is approximately $x'/V$. A prime denotes differentiation with respect to the time, $t$. The responses in the right-hand column show the amplitude of the steady-state oscillations. These periodic motions are stable when represented by a solid line, unstable when represented by a broken line. Hopf bifurcations on the trivial solution are denoted by H, and away from the trivial path stable and unstable oscillatory regimes meet at cyclic folds. Fast dynamic jumps are indicated by vertical arrows.

For case (a) the wavy arrow denotes a slightly turbulent wind (elsewhere the wind is steady). The 2:1 rectangular cross-section exhibits a super-critical Hopf bifurcation at H, with a path of stable limit cycles for higher values of the wind speed. In row (b) the square cross-section in a steady wind exhibits at H a super-critical Hopf bifurcation; and the subsequent limit cycles exhibit two cyclic folds and an associated hysteresis cycle. In row (c) a 2:1 rectangle in steady wind exhibits a sub-critical Hopf bifurcation at H from which a fast dynamic jump would carry the system to a large amplitude stable limit cycle (a periodic attractor). The unstable path from H eventually stabilizes at a cyclic fold, giving an overall (dynamic) response akin to the (static) response of many shell-buckling problems.

In the bottom row, (d), a 1:2 rectangle standing across-wind gives no bifurcation from the trivial solution but large amplitude stable and unstable cycles do exist, separated again by a fold.

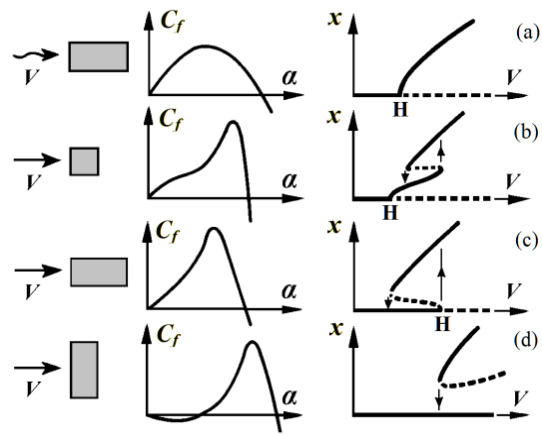

Figure 1 Various aero-dynamic characteristics (first column) and their corresponding dynamic responses (second column) due to Novak [1].

Some of the most familiar examples of galloping arise with engineering cables [6, 7], but we should note that a cable of circular cross-section cannot gallop because the (pure drag) force is in the direction of the resultant wind velocity, and therefore opposes any cable motion. Some cables that can and do gallop are shown in figure 2.

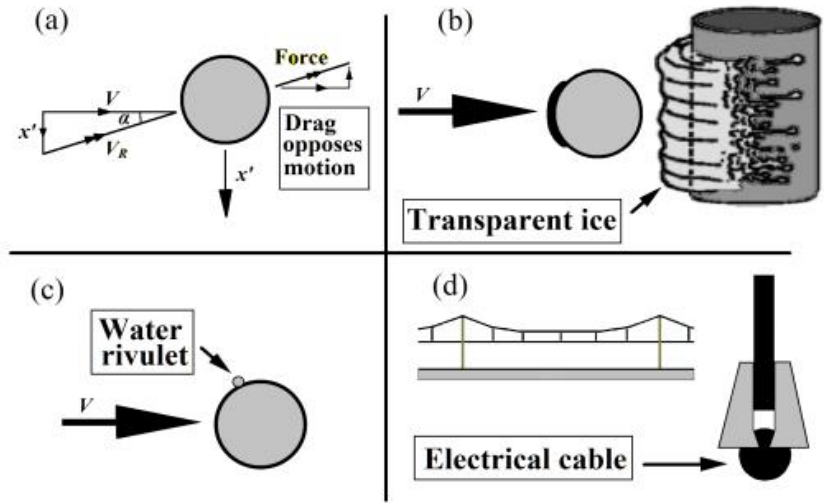

Fig 2. A circular cross-section, (a), cannot gallop because the drag always opposes the motion. Examples that can gallop are: (b) a cable coated with ice; (c) a cable with a rivulet of water, the position of which gives an extra degree of freedom [7]; (d) the electrical wire of an overhead railway catenary, with its notched cross-section.

Galloping problems can also arise in complete structures, such as tower blocks, and here there can be interactions between the wind-induced vibrations and gravity-induced buckling. A classic case was the high-rise Hancock Tower in Boston [8] which had a lot of such problems in its early days. Window panes started falling out, and eventually all 10,344 had to be replaced (the London Shard has 11,000). Occupants suffered from motion sickness, and tuned mass-dampers had to be fitted. There were still problems, however, when a gravitational instability increased the period of vibration from 12 to 16 seconds. The final cure was to add 1,500 tons of diagonal steel bracing, costing $5 million. The tower is still standing today; and still winning architectural prizes for its minimalism!

It is the purpose of this paper to examine the interactions between (Hopf) galloping and (pitch-fork) buckling, remembering that simultaneous failure modes often represent a simplistic, though potentially dangerous, optimal design [9]. We introduce an archetypal model which is non-conservative but autonomous, subjected to time-independent loading by a steadily flowing fluid (air or water). It is designed to exhibit sub-critical bifurcations in both galloping and buckling, both of which will trigger a dynamic jump to a remote stable attractor. When there is more than one candidate attractor, the one onto which the structure settles after the Hopf bifurcation can be *indeterminate* [5, 10]. This is due to the two-dimensional spiralling outset (unstable manifold) of the Hopf, which makes the outcome sensitive to infinitesimally small variations in starting conditions or parameters. This indeterminacy forms the focus of our investigation.

## 2. Archetypal Model for Combined Galloping & Buckling

We consider the archetypal model, shown in figure 3, that we use to study the nonlinear dynamic interactions between galloping and shell-like buckling. A rigid link is pivoted as shown, and held (nominally) vertical by a long spring of stiffness $k$ which is assumed to remain horizontal throughout and is attached to the mass-less rod at a distance $L_2$ from the pivot. We introduce an imperfection into the model by supposing that this spring is initially too short by $y_0$ to hold the unloaded rod exactly vertical. Loaded by the mass $m$ of the grey prism, assumed concentrated at a point on the mass-less rod at a distance $L_1$ from the pivot, this model will exhibit a sub-critical pitch-fork bifurcation. The only interaction with the wind is (considered to be) through the grey prism which has a 2:1 section with the longer edges lying in the direction of the wind. As we have seen, such a prism was analysed by Novak, and shown to exhibit galloping at a sub-critical Hopf bifurcation. The rotational deflection of the rod is written as $x$, and a prime denotes differentiation with respect to the time, $t$.

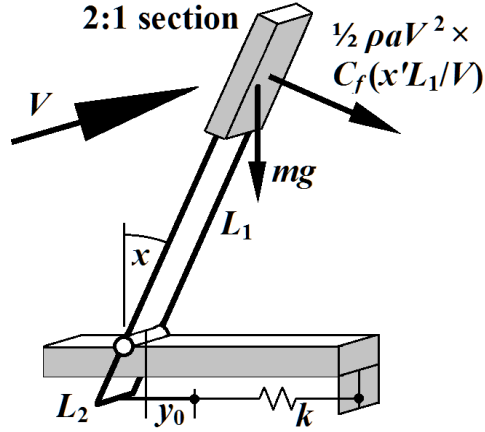

Fig 3. The archetypal model introduced for interaction studies between galloping and shell-like post-buckling. Note that the spring is assumed to be sufficiently long so that it can be assumed to remain horizontal. Meanwhile, as an 'imperfection' in manufacture, the spring is assumed to be initially too short, by $y_0$, to hold the unloaded rod exactly vertical.

The acceleration due to gravity is $g$, the flowing fluid (air, say) has velocity $V$ and density $\rho$ while the frontal area of the prism is written as $a$. The aerodynamic coefficient, $C_f(x'L_1/V)$, is a function of the ratio of the lateral prism speed, $x'L_1$, (assumed to be uniform, and given by the velocity of the centre of gravity) to the wind velocity $V$. This function is typically obtained from wind-tunnel tests in which the *stationary* prism is tilted at a small angle, $\alpha \approx x'L_1/V$, to the direction of flow. Note that in using this *quasi-static* approach pioneered by Novak [1] we are implicitly assuming that the motion of the body is slow compared to the motions of the passing fluid. The equation of motion of the model is

$$mL_1^2 x'' + kL_2(y_0 + L_2 \sin x) \cos x = mgL_1 \sin x + \tfrac{1}{2}\rho a V^2 L_1 C_f(x'L_1/V) \quad (1)$$

We define the following parameters (and hence-forth often omit the word 'parameter'):

*load* parameter, $\Lambda := g/L_1$

*spring* parameter, $B := kL_2^2/mL_1^2$

*imperfection* parameter breaking the pitch-fork symmetry, $e := y_0/L_2$

*velocity* parameter $v := V/L_1$

*pre-factor* for the aerodynamic coefficient, $p := \rho a L_1/m$ (always taken as 0.1)

*forcing* function of the aerodynamics $C_f(x'/v)$

*damping* of the structure, $r$ (always taken as 0.1).

We then have

$$x'' + rx' + B(e + \sin x)\cos x = \Lambda \sin x + \tfrac{1}{2}pv^2 C_f(x'/v) \tag{2}$$

where we have added in the linear damping of magnitude $r$.

### 3. The Pitch-fork bifurcation

Under static conditions (and therefore no aerodynamic forces), with no imperfection and ignoring the trivial solution, we have

$$\Lambda = B \cos x \tag{3}$$

which is the falling sub-critical post-buckling path, $\Lambda(x)$, emerging from the trivial solution at the buckling load, $\Lambda^C = B$. Adding an imperfection we get the well-known imperfection sensitivity [11-14] sketched in figure 4.

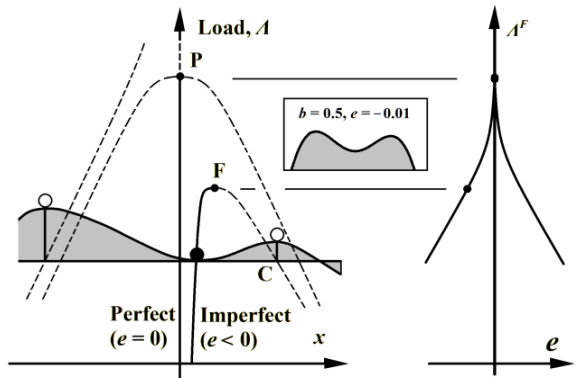

Fig 4. The asymmetric potential well governing the post-buckling behaviour in the presence of a (negative) imperfection.

On the left we have the post-buckling behaviour for a perfect and an imperfect system, with the asymmetric potential well sketched for the latter. The corresponding two-thirds power-law cusp of imperfection sensitivity is shown on the right. The central insert shows the actual shape of the well for the present system at the parameter values of some of our later studies.

Without any loss of generality, we can (by suitable scaling of time) set one of our parameters equal to unity, and for the rest of the paper we set $\Lambda = 1$. So as the pitch-fork parameter we will now use the spring-stiffness $B$ for which $B^C = 1$, and use in particular the *combined stiffness* parameter

$$b := B - 1 = B - B^C \tag{4}$$

as a measure of the 'effective' spring-plus-load stiffness. In the new dimensionless coordinates we have

$$x'' + rx' + (1+b)(e + \sin x)\cos x = \sin x + \tfrac{1}{2}pv^2\, C_f(x'/v)$$

where we keep $p = r = 0.1$ throughout our study. For the forcing function $C_f$ we choose a piecewise polynomial that qualitatively approximates the measurements by Novak [1]:

$$C_f(x') = p(8x')$$

and

$$p(y) = (2/15)y + y^3/3 - y^4/10 - y^5/15 \text{ for } y \geq 0 \text{ and } p(y) = -p(-y) \text{ for } y < 0$$

(see figure 5(a) for the shape of $C_f$). Note that $b = 0$ at the pitch-fork, and positive $b$ measures how far we are away from buckling (on the *stable* trivial path).

## 4. The Hopf bifurcation

Turning to the Hopf bifurcation, we must first think about the aero-dynamic curve, $C_f(x'/v)$, for which we have adopted an analytic function that closely fits Novak's experimentally determined form that we saw in figure 1. This is more suitable for our theoretical work than the power series that Novak used to fit the experiments for $x'/v$ positive, since reflecting this for negative $x'/v$ (as Novak did) gives rise to a singularity at the origin. Our form is shown in figure 5(a).

The Hopf bifurcation arises when the total effective linear damping vanishes, namely when

$$rx' = \tfrac{1}{2}pv^2 x'\, dC_f(x'/v)/dx' \tag{5}$$

Now the derivative of our employed function is $dC_f(x'/v)/d(x'/v) = 1.067$, giving $dC_f(x'/v)/dx' = 1.067/v$, so the value of $v$ at the Hopf bifurcation, $v^H$, is given by

$$v^H = 2r / 1.067\, p \tag{6}$$

With as $r = 0.1$ and $p = 0.1$ (as throughout the paper) we have $v^H =$ **1.875** and the full nonlinear response for these, determined by direct numerical simulation, is displayed in figure 5(b).

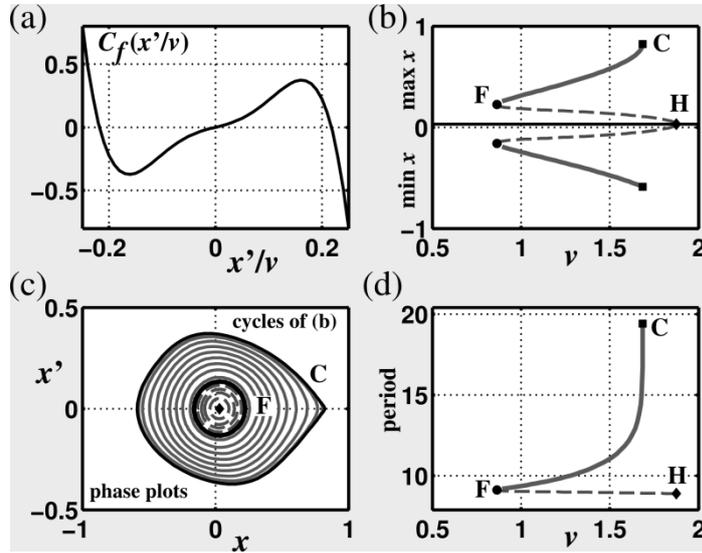

Fig 5 Galloping limit cycles triggered by the Hopf bifurcation at H. A change of stability of the cycles is seen at the cyclic fold, F. Parameter values: *e* = -0.01, *b* = 0.5, $v^H$ = 1.875.

Here we see the sub-critical Hopf bifurcation at H generating the trace of unstable (dashed) cycles which become stable (solid) cycles at the cyclic fold, F. The graph shows the maximum and minimum values of *x(t)* for the steady state galloping oscillations against the wind velocity, *v*. The figure is drawn for an imperfect system with $e = -0.01$, which explains the asymmetry about the *v* axis, and in particular why H does not lie precisely at *x* = 0. Notice that the result $v^H =$ **1.875** is independent of the imperfection. The localized curving of the path of stable cycles at C signifies its approach to the nearer of the two post-buckling equilibrium states (state C of figure 4).

All the cycles (stable and unstable) of figure 5(b) are superimposed on the phase-space portrait of 5(c) where the sharp point of the outer orbit corresponds to the proximity of C. Finally figure 5(d) shows the variation of the periodic times of the cycles traced in 5(b). The period is tending to infinity as the final cycles approach the hill-top equilibrium, C.

## 5. Sequence of phase portraits of the complete model

Looking finally at the complete model, with both wind and gravity loading, we show in figure 6 a sequence of phase portraits for fixed gravity loading and fixed imperfection.

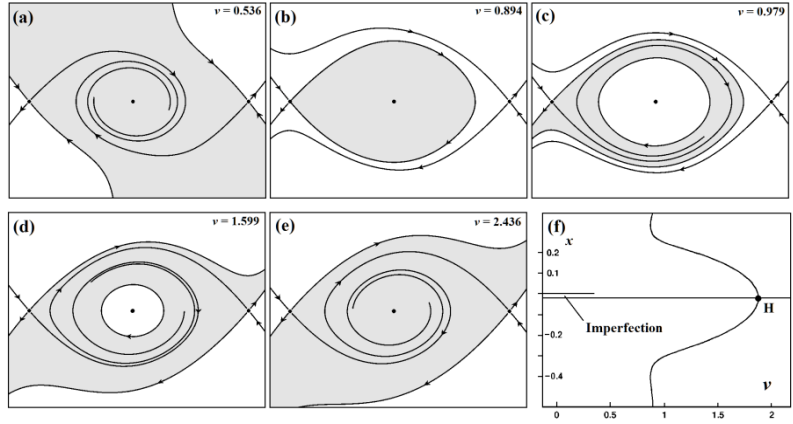

Fig 6. (a) to (e): sequence of phase portraits for an imperfect system below the pitch-fork as the wind speed, *v*, increases. (f): an overview of the corresponding response, showing the maximum and minimum values of the displacement, *x*, against the wind velocity *v*. The Hopf bifurcation is at point H. The fixed parameter values are *e* = 0.003, *b* = 0.175. The grey shading has no technical meaning; it is chosen simply to help the reader see the main features of the diagrams.

In figure 6(a) the topology of the portrait is not yet significantly affected by the wind. Portrait (b) shows a homoclinic connection which together with a very localised fold creates the unstable cycle seen in portrait (c). Notice that disturbances of this cycle generate escape only to the left over the lower potential barrier (corresponding to the *positive* value of *e*). Between portraits (c) and (d) a heteroclinic connection alters the topology, so that in (d) the escape is indeterminate, being either to the left or right. In (e), past the Hopf bifurcation, the central point is unstable and is likewise indeterminate, with disturbances generating escapes over either of the potential hill-tops.

## 6. The co-dimension-two event with symmetry (Kuznetsov)

Before starting our analysis of the co-dimension-three event that governs our symmetry-breaking model, it is useful to look at the unfolding of the *symmetric* case given by Kuznetsov [15]. He takes the normal form of the symmetric Hopf-pitchfork coalescence as

$$x'' - w\,x' - x^2 x' - p\,x - x^3 = 0 \qquad (7)$$

which contains no symmetry-breaking imperfection.

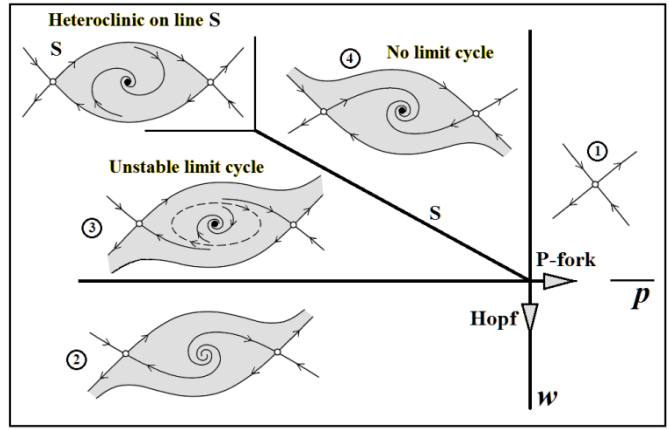

Fig 7. Kuznetsov's complete 2D unfolding of the Hopf-pitchfork, within the family of symmetric systems.

His complete (codimension-2) unfolding of the singularity, in the space of the two control parameters ($w, p$) is shown in figure 7. Notice that the only attractors are the trivial points for $w < 0$. The significant event in this diagram is the saddle connection that occurs on line S, which separates two regions of parameter space, one containing an unstable limit cycle which is destroyed on crossing S.

## 7. Co-dimension-three event of our model

Guided by this 2D unfolding of the symmetric case, we now proceed to fully unfold the compound singularity exhibited by our model in the 3D parameter space of our stiffness parameter, *b* (an inverse representation of the pitch-fork loading), our wind velocity, *v*, and our symmetry-breaking parameter, *e*. This compound bifurcation has been called the Takens-Bodganov Cusp [15, 16]. Note that these authors study the unfolding of the centre-saddle-centre case as opposed to our saddle-centre-saddle case illustrated in fig 7.

The result is shown in figure 8, which gives two views of the same ellipsoid in parameter space.

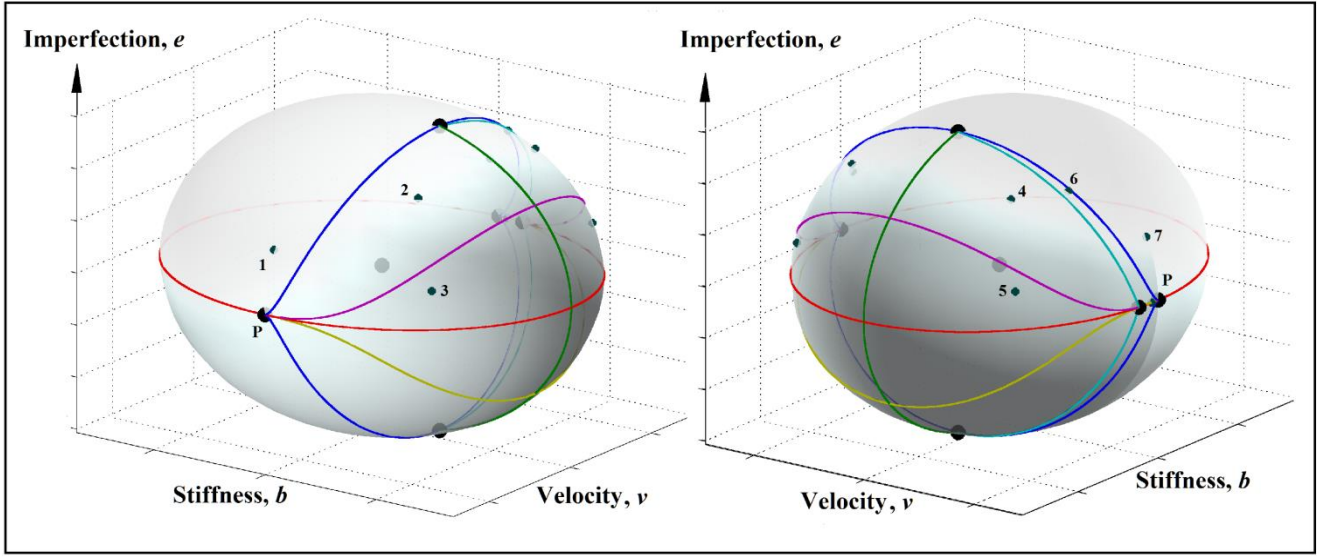

Fig 8. Intersections of the 2D bifurcation surfaces with an ellipsoid in the 3D parameter space. Colour will be available on-line.

The figure shows an ellipsoid surrounding the critical point in the space spanned by the stiffness parameter, *b*, the velocity parameter, *v*, and the imperfection parameter, *e*. It is defined by the polar coordinates (*ϕ*, *ψ*) according to the transformation

$$\begin{aligned} v &= v^H + R^2 \cos \pi\phi \cos \psi\pi/2 \\ b &= b_0 + R \sin \pi\phi \cos \psi\pi/2 \\ e &= e_0 + R^3 \sin \psi\pi/2 \end{aligned} \quad (8)$$

with $v^H$=1.875 (as derived in equation (6)), $b_0 = 0$ and $e_0 = 0$.

Here the 'small' radius parameter, *R* is taken nominally as 0.2, but for clarity the picture is not to scale. The resulting image does not change qualitatively for smaller or greater *R*. The coloured arcs drawn on the ellipsoid show its intersection with the various bifurcation surfaces emerging from the origin, which are better understood in the unfolded ellipsoidal surface of figure 9 when projected into the parameter plane (*ϕ*, *ψ*).

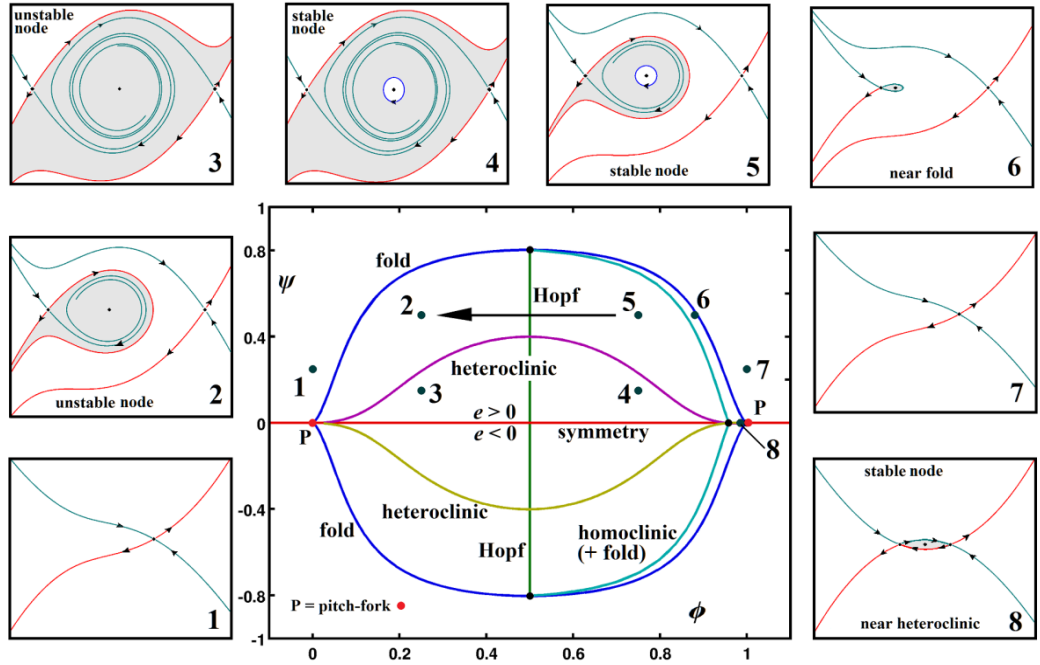

Fig 9. Unfolding of the compound (codimension-three) galloping-buckling bifurcation. A projected view of the ellipsoidal surface showing the bifurcation arcs of figure 8. Colour will be available on-line.

In this view, the pitch-fork bifurcation appears twice, at P, where the (blue) arc of static folds exhibits a very localized cusp. The Hopf bifurcation (dark green) occurs on the vertical axis. There are now two types of saddle connection, a homoclinic (with itself) and a heteroclinic (between two distinct unstable saddle equilibria). Crossing the heteroclinic (purple) arc takes us, for example, from portrait (2) where the galloping system escapes only to the left to portrait (3) where the outcome is indeterminate, depending sensitively on the starting condition near the node. Remember, here, that the asymmetry, $e$, varies as we move over the ellipsoidal surface; it is positive above the red symmetry line, making the lower escape barrier on the left, unlike as in figure 4 where $e$ was negative. We explore this sensitivity more fully in the following section. Meanwhile, crossing the homoclinic (light green) arc transforms portrait (5) into portrait (6). Notice that portrait (6) is very close to the fold (blue) arc, crossing which gives a portrait such as (7) with only one equilibrium fixed point. A feature not visible in figure 9 (but present) for small radius $R$ is the fold of limit cycles (very close to the homoclinic). This fold of limit cycles must exist because the periodic orbit born in the homoclinic is asymptotically stable. The fold of limit cycles is more clearly visible in figures 5(b, d) and 6(f). Correspondingly, there should be a phase portrait in figure 9 between portraits 9(5) and 9(6) with two coexisting limit cycles (the inner one unstable, the outer one stable); this could not be conveniently shown in figure 9.

## 8. Ramped velocity and indeterminate outcomes

Our results for ramping the velocity $v$ as a linear function of time according to the equation

$$v(t) = v(0) + \gamma\, t \tag{9}$$

are shown in figure 10. The approximate right-left symmetries of these graphs about the Hopf line (particularly pronounced in the pictures of the second column) are discussed briefly in the Appendix.

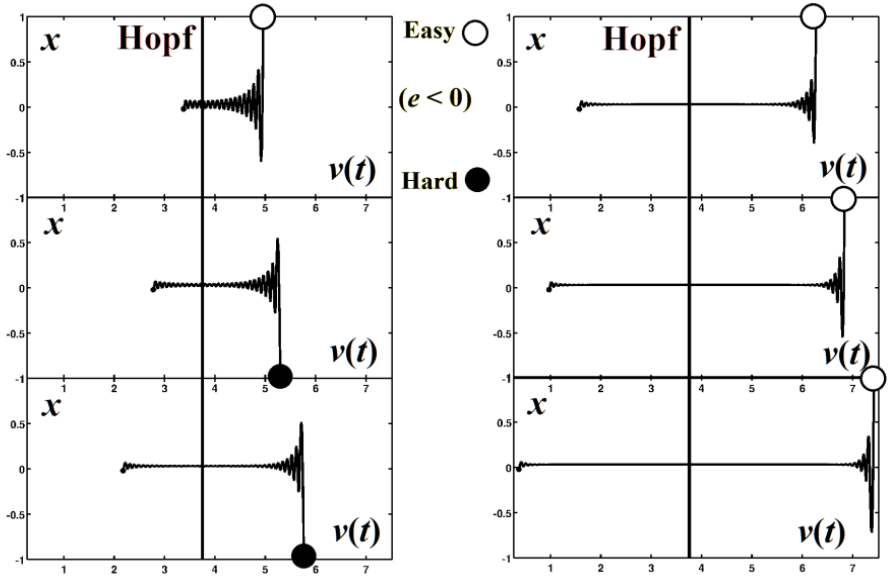

Fig 10. Time-series results with a ramped velocity, showing tunnelling and indeterminacy. The velocity, $v$, is ramped from varied distances below the Hopf bifurcation. The following are held constant throughout the figure: stiffness, $b = 0.5$, imperfection $e = -0.01$, $x(0) = x_{eq} - 0.05$, $x'(0) = 0$, $\gamma = 0.01$, where $x_{eq}$ is the central equilibrium value of $x$.

This study is for negative $e$, making it 'easier' for the system to escape to positive large $x$, but the parameters are such, as in portrait (3) of figure 9, that escape from the un-ramped Hopf bifurcation is indeterminate being possible in either direction (towards $x$ positive or $x$ negative).

We notice first the considerable 'tunnelling' through the Hopf bifurcation which arises because the small disturbance from equilibrium takes time to grow under the light negative effective damping just after the steady-state Hopf velocity. This tunnelling increases as the runs are started earlier and earlier, because the longer time interval under positive damping ensures that $x$ and $x'$ have decreased closer and closer towards the origin before $v$ reaches $v^H$.

Next we observe that some runs escape over the lower hill-top equilibrium ($x > 0$) while others escape over the higher hill-top ($x < 0$). The relative hill-height for this value of $b$ is shown as an insert in figure 4. A further study of this indeterminacy is shown in figure 11.

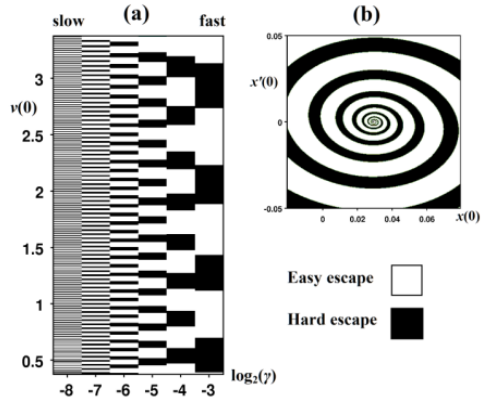

Fig 11. Illustration of indeterminate outcomes for parameter values $b = 0.5$, $e = -0.01$, $r = 0.2$, $v^H = 3.75$.

In figure 11(a) we display the outcomes, in terms of easy escape over the lower barrier (white) or hard over the higher barrier (black) resulting from different values of $v(0)$ for different ramping rates corresponding to the six integer values of $\log_2(\gamma)$. The fixed starting conditions in (a) are $x(0) = x_{eq} - 0.05$, $x'(0) = 0$ where $x_{eq}$ is the equilibrium value of $x$. In figure 11(b) we show, again in black or white, the outcomes in the space of $[x(0), x'(0)]$ for fixed values of $v(0) = v^H/2$ at $\gamma = 0.01$.

While our numerical simulations of parameter-ramping through a Hopf bifurcation are adequate for the case when the location of the equilibrium does not depend on the drifting parameter (which is the case in our present investigation), we should note that the general case is more subtle. In particular, the value of the drifting parameter at which the trajectory starts to grow noticeably exponentially depends not only on the starting parameter ($v(0)$ in our case) but also on properties of the right-hand side of the equation. See [17] for a mathematical treatment, and [18] for some of the typical observations.

## 9. Concluding remarks

We have proposed and studied an archetypal model to explore the nonlinear dynamic interactions between galloping at an incipient sub-critical Hopf bifurcation of a structure with shell-like buckling behaviour. Optimal designs often call for a simultaneity of failure modes, but nonlinear interactions can then be dangerous [9]. The compound bifurcation corresponding to simultaneous galloping and buckling is the so-called Takens-Bodganov Cusp, and we have made a full unfolding of this codimension-3 bifurcation for the model to explore the adjacent phase-space topologies.

The indeterminacy of the outcome, that we find for both quasi-static and ramped loadings, should certainly be noted by design engineers. It will be interesting to see if the various approaches of analysis and control of safe basins of attraction [19, 20] pioneered by Giuseppe Rega and Stefano Lenci and can play a role in interactions of the present type.

## Appendix

In the second column of figure 10 we have noted the approximate right-left symmetries of the time-series about the Hopf line. A simple approximate analysis of this, valid in our case because change of $v$ does not change the position of the equilibrium solution, can be written down as follows.

Let $c(v)$ be the real part of the complex eigenvalue of the equilibrium for a given $v$. For $v < v^H$ we have $c(v) < 0$, while for $v > v^H$ we have $c(v) > 0$. If $v$ drifts slowly then the amplitude of a disturbed trajectory will initially decay exponentially, with rate $c(v)$, and then grow.

The picture of the time-series will be symmetric if the real part $c(v)$ is odd about $v^H$, and we remember that close to a Hopf bifurcation $c(v)$ increases linearly. So, close to the Hopf bifurcation the time series is always approximately symmetric, as long as the displacements from the equilibrium state are small. In general, the approximation formula is

$$d(v) = d_0 \; \exp \left[ \int_{v(0)}^{v} c(w) \; dw / \gamma \right]$$

where, ignoring rotation, $d(v)$ is the amplitude of the displacement when the wind speed is $v$. The drift speed is $\gamma$, the initial wind speed is $v_0$ (written as $v(0)$ in the equation), and $d_0$ is the value of $d$ at $v = v_0$. Note that the exponential can become extremely small, so noise and round-off determine the point where growth becomes noticeable in the figure.

## Acknowledgement

The research of J.S. is supported by EPSRC Grant EP/J010820/1.

x Revision FFCL x